\documentclass[12pt,english]{article}
\usepackage[T1]{fontenc}
\usepackage[latin9]{inputenc}
\usepackage{geometry}
\usepackage{babel}
\usepackage{float}
\usepackage{amsmath}
\usepackage{amssymb}
\usepackage{graphicx}
\usepackage{hyperref}
\usepackage{breakurl}

\makeatletter

\providecommand{\tabularnewline}{\\}

\usepackage{babel}

\makeatother

\begin{document}

\title{Quantized-CP Approximation and Sparse Tensor Interpolation of Function Generated Data}

\author{Boris N. Khoromskij\thanks{Max-Planck-Institute for Mathematics in the Sciences, 
Inselstr. 22-26, 04103 Leipzig, Germany (bokh@mis.mpg.de).}, ~
Kishore K. Naraparaju\thanks{Birla Institute of Technology and Science Pilani, Hyderabad Campus,
Hyderabad, India (thimmaki@gmail.com)} 
\thanks{This work is partially supported by National Board of Higher Mathematics, DAE, India.}
~~and Jan Schneider\thanks{Wests\"achsische Hochschule Zwickau, Dr. Friedrichs-Rings 2A, Zwickau,
Germany (jan.schneider@fh-zwickau.de). }}

\date{~}
\maketitle
\begin{abstract}
In this article we consider the iterative  schemes to compute the 
canonical (CP) approximation of quantized
data generated by a function discretized on a large uniform grid in an 
interval on the real line. This paper continues the research on
the QTT method   \cite{KhQuant:09} developed for the tensor train (TT) approximation 
of the quantized images of function related data. 
In the QTT approach the target vector of length $2^{L}$ is reshaped to a $L^{th}$
order tensor with two entries in each mode (Quantized representation)
and then approximated by the QTT tenor including $2r^2 L$ parameters, where $r$ is the maximal TT rank. 
In what follows, we consider the Alternating Least-Squares (ALS) iterative scheme to compute
the rank-$r$ CP approximation of the quantized vectors, 
which requires only $2 r L\ll 2^L$ parameters for storage.
In the earlier papers \cite{KhorSurv:10} such a representation was 
called Q$_{Can}$ format, while in this paper we abbreviate it as the QCP representation. 
 We test the  ALS algorithm to calculate the QCP approximation on various functions,
and in all cases we observed the exponential error decay in the QCP rank. 
The main idea  for recovering a discretized function
in the rank-$r$ QCP format using the reduced number the functional samples, 
calculated only at $O(2rL)$ grid points, is presented. 
The special version of ALS scheme for solving the arising minimization problem
is described. This approach can be viewed as the sparse QCP-interpolation method
that allows to recover all $2r L$ representation parameters of the rank-$r$ QCP tensor.
Numerical examples show the efficiency of the QCP-ALS type iteration and indicate 
the exponential convergence rate in $r$. 
\end{abstract}


\section{Introduction}\label{sec:Introd}

In many applications, the approximation or integration of functions 
inheriting the properties of $e^{-kx},\,e^{-kx^{2}},\,e^{-k|x|},\, \sin(kx)$ 
or $1/|x|^\alpha$ on an interval in $\mathbb{R}$ as well as functions depending on many parameters
leads to the challenging numerical problems. Often a very
fine grid is required to approximate sharp functions like the
Gaussians $e^{-kx^{2}}$ for large values of $k$, highly oscillating
functions or functions with multiple local singularities or cusps arising, for example,
as the solution of PDEs discretized on fine spatial grid. 
The storage of the function values as well as simple
arithmetic's operations on data arising from sampling on large grids may easily 
become non-tractable. The additional difficulty arises if each function evaluation 
has very high cost, say, related to the solution of large linear system 
or spectral problem as well as to solving complicated PDE.

The quantics-TT (QTT) tensor approximation method,  introduced and analyzed 
 in \cite{KhQuant:09} for some classes of discretized functions, is now 
a well established technique for data compression of long function generated
vectors. It is based on the low-rank tensor approximation to the quantized image
of a vector, where the tensor train (TT) format \cite{osel} was applied 
to the quantized multi-fold image.
Based on the quantization (reshaping) of a long $2^L$-vector to 
a $L^{th}$-order tensor (Quantics)
the consequent QTT tensor approximation has been proven to have low 
TT rank for a wide class of functional vectors.
We refer to \cite{Osel-TT-LOG:09} where the TT approximation to the reshaped $2^L\times 2^L$
Laplacian type matrices was considered and analyzed numerically. 
The QTT tensor parametrization requires $O(2r^2 L)$
storage size where $r$ is the upper bound on the TT rank parameters.
Some examples on the successful application of the QTT tensor approximation to the solution 
of PDEs and in stochastic modeling can be found in 
\cite{BeDoKh2_BSE2:16,DoKh_QTTT:12,KaOsRaSch:16,KaRaSch:13,KhSautVeit:11} and in 
\cite{Boris2,boris3,KhSchw1:09,BokhSRep:15,VeBoKh:Ewald:14,DolKhLitMat:14,VeKhSchn_QTTHF:11},
among others.

The present paper  continues the research on
the QTT tensor approximation method  \cite{KhQuant:09} based on the use of TT format.
In what follows we investigate the numerical schemes to compute canonical (CP) tensor approximation
of the quantized tensor. This data format was introduced in \cite{KhorSurv:10} under the name
Q$_{Can}$  tensor representation. In this paper,
we shall abbreviate the notion Q$_{Can}$ as the QCP format.
First, we briefly recall the main construction along the line of the QTT approximation.
A given
vector of size $N=2^{L}$ is reshaped (quantized) by successive dyadic folding
to a $\underset{L}{\underbrace{2\times2\times\ldots\times2}}$ array.
The rank $r$ representation of this tensor in the canonical format
reduces the number of representation parameters from $2^L$ down to $2 r L$,
 which is smaller than for the QTT format, characterized by the storage size $O(2r^2 L)$.
The following simple example shows why the QCP approximation of a
vector does a job by reducing the number of representation parameters.

Let $f(x)=e^{-x}$ in $[0,1]$. Consider the nodes $0,h,2h,\ldots,15h$
on the interval $[0,1],$ where $h$ is the step size of the uniform
grid. The function values at these discrete points form a vector
${\bf q}=[1,q,q^{2},\ldots,q^{15}]$ of length $N=2^{4}$, where $q=e^{-h}.$
Now reshape the vector ${\bf q}$ to obtain a $4$th order tensor 
${\bf Q}\in \mathbb{R}^{2\times 2\times 2\times 2}$ that is the 
quantized image of ${\bf q}$.
Following \cite{KhQuant:09}, we recall that the CP rank of
this tensor is 1, and the corresponding explicit CP tensor representation reads as 
\begin{eqnarray*}
{\bf Q}=\left(\begin{array}{c}
1\\
q
\end{array}\right)\otimes\left(\begin{array}{c}
1\\
q^{2}
\end{array}\right)\otimes\left(\begin{array}{c}
1\\
q^{4}
\end{array}\right)\otimes\left(\begin{array}{c}
1\\
q^{8}
\end{array}\right).
\end{eqnarray*}

One can see that the whole vector of size $N=2^4$ is represented only by 4 parameters,
that means the logarithmic complexity scaling $\log_{2}N.$ In general, the CP rank of the 
quantized image of a vector
of length $2^{L}$ generated by $f(x)=e^{-\lambda x}$ is $1$ and it is
represented by only $\log_{2}2^{L}=L$ parameters, such that its explicit one-term representation
is given by \cite{KhQuant:09} 
\begin{eqnarray*}
{\bf Q}=
\bigotimes_{p=1}^{L}\left(\begin{array}{c}
1\\
q^{2^{p-1}}
\end{array}\right).
\end{eqnarray*}
We say that the CP rank of the quantized image of the discretized function $e^{-\lambda x}$
is $1.$ In general, the exact CP rank of rather simple functions 
like $e^{-x^{2}},\sin\,k\pi x,|x|^\alpha,...$ etc. is not known. 
Construction and complexity analysis of ALS-type algorithms for computing the QCP approximation of 
quantized functions is the main aim of this article.

It is worth to note that the rank-$r$ QCP tensor is represented only by small number of
 parameters, $2Lr$, whereas the QTT format based on rank $r$ TT tensors 
is parametrized by $O(2Lr^{2})$ numbers as it was already mentioned.
Based on this observation, we propose the QCP interpolation scheme using only small number
of functional calls (of the order of $O(2Lr)$) which recovers the quantized tensor image.
This concept leads to the promising enhancement of the QPC 
approximation of the complete $2^L$-vector because the small number of parameters
in the arising minimization problem. For the practical implementation, 
we introduce the ALS type scheme to compute the sparse QCP interpolant.


In numerical examples, we test the ALS iterative scheme implementing 
the CP approximation on $15$th order tensors representing $2^{15}$-vectors generated by
various functions. In all cases we observe the exponentially fast error decay
in the QCP rank.
Notice that the traditional ALS algorithm for CP tensors and its enhanced versions have 
been discussed in many articles 
\cite{carrollchang,common,domanov,espighack,golubvanloan,harshman,tamarakolda,KhSchw1:09,usch,vanloan,RoHoSchn_ALS:10}.
Regularized ALS scheme was considered in \cite{likinder,NVASCALAT}.

The efficient representation and multilinear algebra of large multidimensional
vectors (tensors) can be based on their low-rank tensor approximation
by using different tensor formats. We refer to reviews on the multilinear algebra
\cite{tamarakolda,Hack_Book12,lars,Ours,OsSaTy_LATen:09} 
and to recent surveys on tensor numerical methods and their application in scientific
computing \cite{KhorSurv:10,KhorSurv:14,VeKhorTromsoe:15,CichOseletc:2016}.

The rest of the paper is organized as follows. Some auxiliary technical 
results concerning the ALS-canonical algorithm are presented in section 2. 
Section 3 describes the particular QCP-ALS scheme 
which uses the complete information about the tensor. 
In section 4, we calculate the QCP approximation of some
selected functions which appears in various applications. The main idea and basic ALS scheme for
computing the QCP approximation by using the information on
only few entries of the target vector is described in Section 5, where 
the numerical illustrations are also presented.
The approximation by using incomplete data can be viewed as the sparse QCP interpolation 
of function generated vectors.
Some useful notations, definitions and a simple example of the scheme for QCP
approximation of a $4^{th}$ order tensor are given in Appendix.

In this article, we often use MATLAB notations, for example 
${\bf X}=\mbox{reshape}({\bf x},2,...,2)$. 
The Frobenius norm of a tensor
${\bf X}=[x_{i_{1}i_{2}\ldots i_{d}}]\in\mathbb{R}^{n_{1}\times n_{2}\times\ldots\times n_{d}}$
is defined as the square root of the sum of squares of all its elements $x_{i_{1}i_{2}\ldots i_{d}}$,
i.e, 
\begin{eqnarray*}
\left\Vert {\bf X} \right\Vert _{F}=
\sqrt{\sum_{i_{1}=1}^{n_{1}}\sum_{i_{2}=1}^{n_{2}}\ldots
\sum_{i_{d}=1}^{n_{d}}x_{i_{1}i_{2}\ldots i_{d}}^{2}}.
\end{eqnarray*}

\section{Technical results for ALS}

For convenience and better understanding of our notation and some
technical results in the upcoming sections, we provide the proofs of some basic
results (see \cite{vanloan}), which are useful in the construction of ALS algorithm, 
and give some simple examples of ALS iteration for canonical approximation of a tensor. 

A complete description of each iteration step of ALS for a
fourth order tensor $\chi$ is given in the appendix A3. In the steps 1,~2,~3
and 4 of A3, we need to obtain products like $(C\odot B\odot A)^{T}(C\odot B\odot A)$
and $(C\odot B\odot A)^{T}\chi(:,:,:,1)$, where $\odot$ is the Khatri-Rao
product of matrices defined in appendix A2. An efficient
way of obtaining these products is described below. First, we show it
for products which appear in the canonical approximation of a
third order tensor and then generalize it for $d^{th}$ order tensors.

\subsection{Fast evaluation of $\left(B\odot C\right)^{T}\left(B\odot C\right)$}

In the Alternating Least Squares method, to get the Canonical approximation
to a third order tensor, we need to compute $\left(B\odot C\right)^{T}\left(B\odot C\right).$
This usually requires $O(r^{2}n_{1}n_{2})+n_{1}n_{2}r$ arithmetic operations.
Now we show the efficient way to compute it.

\subsubsection*{Lemma 1.\textmd{ If 
$B=[\mathbf{b}_{1}|\mathbf{b}_{2}|\ldots|\mathbf{b}_{r}]\in\mathbb{R}^{n_{1}\times r}$
and $C=[\mathbf{c}_{1}|\mathbf{c}_{2}|\ldots|\mathbf{c}_{r}]\in\mathbb{R}^{n_{2}\times r},$
then 
\begin{eqnarray*}
\left(B\odot C\right)^{T}\left(B\odot C\right)=\left(B^{T}B\right)\circ\left(C^{T}C\right)
\end{eqnarray*}
where $\circ$ denotes the Hadamard product of matrices. }}

\textbf{Proof:} As defined in appendix A2,
\begin{eqnarray*}
B\odot C=[\mathbf{b}_{1}\otimes\mathbf{c}_{1}|\mathbf{b}_{2}\otimes\mathbf{c}_{2}|\ldots|\mathbf{b}_{r}
\otimes\mathbf{c}_{r}]\in\mathbb{R}^{n_{1}n_{2}\times r}.
\end{eqnarray*}
So 
\begin{eqnarray*}
(B\odot C)^{T}(B\odot C)=\left[\begin{array}{c}
(\mathbf{b}_{1}\otimes\mathbf{c}_{1})^{T}\\
(\mathbf{b}_{2}\otimes\mathbf{c}_{2})^{T}\\
\vdots\\
(\mathbf{b}_{r}\otimes\mathbf{c}_{r})^{T}
\end{array}\right]_{r\times n_{1}n_{2}}[\mathbf{b}_{1}\otimes\mathbf{c}_{1}|\mathbf{b}_{2}\otimes\mathbf{c}_{2}|
\ldots|\mathbf{b}_{r}\otimes\mathbf{c}_{r}]_{n_{1}n_{2}\times r}
\end{eqnarray*}
\begin{eqnarray*}
=\left[\begin{array}{ccccc}
(\mathbf{b}_{1}\otimes\mathbf{c}_{1})^{T}(\mathbf{b}_{1}\otimes\mathbf{c}_{1})\,\, & (\mathbf{b}_{1}\otimes
\mathbf{c}_{1})^{T}(\mathbf{b}_{2}\otimes\mathbf{c}_{2}) & ... & ... & (\mathbf{b}_{1}\otimes\mathbf{c}_{1})^{T}(\mathbf{b}_{r}\otimes\mathbf{c}_{r})\\
(\mathbf{b}_{2}\otimes\mathbf{c}_{2})^{T}(\mathbf{b}_{1}\otimes\mathbf{c}_{1})\,\, & (\mathbf{b}_{2}\otimes
\mathbf{c}_{2})^{T}(\mathbf{b}_{2}\otimes\mathbf{c}_{2}) & ... & ... & (\mathbf{b}_{2}\otimes\mathbf{c}_{2})^{T}
(\mathbf{b}_{r}\otimes\mathbf{c}_{r})\\
\vdots & \vdots &  &  & \vdots\\
(\mathbf{b}_{r}\otimes\mathbf{c}_{r})^{T}(\mathbf{b}_{1}\otimes\mathbf{c}_{1})\,\, & (\mathbf{b}_{r}\otimes
\mathbf{c}_{r})^{T}(\mathbf{b}_{2}\otimes\mathbf{c}_{2}) & ... & ... & (\mathbf{b}_{r}\otimes\mathbf{c}_{r})^{T}(\mathbf{b}_{r}\otimes\mathbf{c}_{r})
\end{array}\right]\\
\\
=\left[\begin{array}{ccccc}
(\mathbf{b}_{1}^{T}\otimes\mathbf{c}_{1}^{T})(\mathbf{b}_{1}\otimes\mathbf{c}_{1})\,\, & (\mathbf{b}_{1}^{T}\otimes
\mathbf{c}_{1}^{T})(\mathbf{b}_{2}\otimes\mathbf{c}_{2}) & ... & ... & (\mathbf{b}_{1}^{T}\otimes\mathbf{c}_{1}^{T})(\mathbf{b}_{r}\otimes\mathbf{c}_{r})\\
(\mathbf{b}_{2}^{T}\otimes\mathbf{c}_{2}^{T})(\mathbf{b}_{1}\otimes\mathbf{c}_{1})\,\, & (\mathbf{b}_{2}^{T}\otimes
\mathbf{c}_{2}^{T})(\mathbf{b}_{2}\otimes\mathbf{c}_{2}) & ... & ... & (\mathbf{b}_{2}^{T}\otimes\mathbf{c}_{2}^{T})(\mathbf{b}_{r}\otimes\mathbf{c}_{r})\\
\vdots & \vdots &  &  & \vdots\\
(\mathbf{b}_{r}^{T}\otimes\mathbf{c}_{r}^{T})(\mathbf{b}_{1}\otimes\mathbf{c}_{1})\,\, & (\mathbf{b}_{r}^{T}\otimes
\mathbf{c}_{r}^{T})(\mathbf{b}_{2}\otimes\mathbf{c}_{2}) & ... & ... & (\mathbf{b}_{r}^{T}\otimes\mathbf{c}_{r}^{T})(\mathbf{b}_{r}\otimes\mathbf{c}_{r})
\end{array}\right]\,\,\,\, & \,\text{(see\,\,P1\,\,in\,\,A1)}\\
\\
=\left[\begin{array}{ccccc}
(\mathbf{b}_{1}^{T}\mathbf{b}_{1}\otimes\mathbf{c}_{1}^{T}\mathbf{c}_{1})\,\, & (\mathbf{b}_{1}^{T}\mathbf{b}_{2}\otimes
\mathbf{c}_{1}^{T}\mathbf{c}_{2}) & ... & ... & (\mathbf{b}_{1}^{T}\mathbf{b}_{r}\otimes\mathbf{c}_{1}^{T}\mathbf{c}_{r})\\
(\mathbf{b}_{2}^{T}\mathbf{b}_{1}\otimes\mathbf{c}_{2}^{T}\mathbf{c}_{1})\,\, & (\mathbf{b}_{2}^{T}\mathbf{b}_{2}\otimes
\mathbf{c}_{2}^{T}\mathbf{c}_{2}) & ... & ... & (\mathbf{b}_{2}^{T}\mathbf{b}_{r}\otimes\mathbf{c}_{2}^{T}\mathbf{c}_{r})\\
\vdots & \vdots &  &  & \vdots\\
(\mathbf{b}_{r}^{T}\mathbf{b}_{1}\otimes\mathbf{c}_{r}^{T}\mathbf{c}_{1})\,\, & (\mathbf{b}_{r}^{T}\mathbf{b}_{2}\otimes
\mathbf{c}_{r}^{T}\mathbf{c}_{2}) & ... & ... & (\mathbf{b}_{r}^{T}\mathbf{b}_{r}\otimes\mathbf{c}_{r}^{T}\mathbf{c}_{r})
\end{array}\right].\,\,\,\,\,\,\,\,\,\,\,\,\,\,\,\,\,\,\,\,\,\,\,\,\,\,\,\,\,\,\,\,\,\,\,\,\,\, & \text{(see\,\,P3\,\,in\,\,A1)}
\end{eqnarray*}
Since $\mathbf{b}_{i}^{T}\mathbf{b}_{j}$ and $\mathbf{c}_{i}^{T}\mathbf{c}_{j}$
are scalars, $(\mathbf{b}_{i}^{T}\mathbf{b}_{j})\otimes(\mathbf{c}_{i}^{T}\mathbf{c}_{j})=(\mathbf{b}_{i}^{T}\mathbf{b}_{j})
\cdot(\mathbf{c}_{i}^{T}\mathbf{c}_{j})$.
Therefore, 
\begin{eqnarray*}
(B\odot C)^{T}(B\odot C)=\left[\begin{array}{cccc}
\mathbf{b}_{1}^{T}\mathbf{b}_{1} & \mathbf{b}_{1}^{T}\mathbf{b}_{2} & ... & \mathbf{b}_{1}^{T}\mathbf{b}_{r}\\
\mathbf{b}_{2}^{T}\mathbf{b}_{1} & \mathbf{b}_{2}^{T}\mathbf{b}_{2} & ... & \mathbf{b}_{2}^{T}\mathbf{b}_{r}\\
\\
\mathbf{b}_{r}^{T}\mathbf{b}_{1} & \mathbf{b}_{r}^{T}\mathbf{b}_{2} & ... & \mathbf{b}_{r}^{T}\mathbf{b}_{r}
\end{array}\right]\circ\left[\begin{array}{cccc}
\mathbf{c}_{1}^{T}\mathbf{c}_{1} & \mathbf{c}_{1}^{T}\mathbf{c}_{2} & ... & \mathbf{c}_{1}^{T}\mathbf{c}_{r}\\
\mathbf{c}_{2}^{T}\mathbf{c}_{1} & \mathbf{c}_{2}^{T}\mathbf{c}_{2} & ... & \mathbf{c}_{2}^{T}\mathbf{c}_{r}\\
\\
\mathbf{c}_{r}^{T}\mathbf{c}_{1} & \mathbf{c}_{r}^{T}\mathbf{c}_{2} & ... & \mathbf{c}_{r}^{T}\mathbf{c}_{r}
\end{array}\right]=(B^{T}B)\circ(C^{T}C).\,\,\,\,\,\,\,\,
\end{eqnarray*}
So one can easily show that $B^{T}B$ requires $O(n_{1}r^{2})$ and
$C^{T}C$ requires $O(n_{2}r^{2})$ arithmetic operations. Therefore, the
computational complexity to compute $(B\odot C)^{T}(B\odot C)$ is
$O((n_{1}+n_{2})r^{2})+r^{2}.$

\subsubsection*{Generalization of Lemma 1}

Here we generalize Lemma 1 to more than two matrices.
Let us consider $A_{1},A_{2},\ldots,A_{L}$ to be matrices of the same size
$n\times r.$ Then by recursion one can easily prove that 
\begin{eqnarray*}
\left[A_{L}\odot..\odot A_{i+1}\odot A_{i-1}\odot A_{i-2}\odot..\odot A_{1}\right]^{T}\left[A_{L}\odot..\odot A_{i+1}\odot A_{i-1}
\odot A_{i-2}\odot..\odot A_{1}\right]\,\,\,\,\,\,\,\,\,\,\,\,\,\,\,\,\,\,\,\,\,\,\,\,\,\,\,\,\,\,\,\,\,\,\,\,\\
=\left[\left(A_{L}\odot..\odot A_{i+1}\odot A_{i-1}\odot A_{i-2}\odot..\odot A_{2}\right)\odot A_{1}\right]^{T}\left[\left(A_{L}
\odot..\odot A_{i+1}\odot A_{i-1}\odot A_{i-2}\odot..\odot A_{2}\right)\odot A_{1}\right]\\
=\left(A_{L}\odot..\odot A_{i+1}\odot A_{i-1}\odot A_{i-2}\odot..\odot A_{2}\right)^{T}\left(A_{L}\odot..\odot A_{i+1}\odot A_{i-1}
\odot A_{i-2}\odot..\odot A_{2}\right)\circ A_{1}^{T}A_{1}\,\,\,\,\,\,\,\,\,\,\,\,\,\\
\vdots\,\,\,\,\,\,\,\,\,\,\,\,\,\,\,\,\,\,\,\,\,\,\,\,\,\,\,\,\,\,\,\,\,\,\,\,\,\,\,\,\,\,\,\,\,\,\,\,\,\,\,\,\,\,\,\,\,\,\,\,\,\,
\,\,\,\,\,\,\,\,\,\,\,\,\,\,\,\,\,\,\,\,\,\,\,\,\,\\
=\left(A_{L}^{T}A_{L}\right)\circ\left(A_{L-1}^{T}A_{L-1}\right)\circ\ldots\circ\left(A_{i+1}^{T}A_{i+1}\right)
\circ\left(A_{i-1}^{T}A_{i-1}\right)\circ\ldots\circ\left(A_{2}^{T}A_{2}\right)\circ\left(A_{1}^{T}A_{1}\right).\,\,\,\,\,\,
\,\,\,\,\,\,\,\,\,\,\,\,\,\,\,\,\,\,
\end{eqnarray*}

The computational complexity to compute the above is $O((L-1)nr^{2})$
whereas the direct computation of this product requires $n^{L-1}r+O(r^{2}n^{L-1}).$
So this is much faster. 

\subsection{Fast evaluation of $\left(B\odot C\right)^{T}\,\mathbf{x}$}

In ALS, we also need to compute $\left(B\odot C\right)^{T}\mathbf{x}$
for $\mathbf{x}\in\mathbb{R}^{n_{1}n_{2}}.$ It would require $3n_{1}n_{2}r-r$
arithmetic operations including $n_{1}n_{2}r$ operations for computing $\left(B\odot C\right).$ 
This complexity can be further improved in the following way. \\
 
 \textbf{Lemma 2: }If 
 $B=[\mathbf{b}_{1}|\mathbf{b}_{2}|\ldots|\mathbf{b}_{r}]\in\mathbb{R}^{n_{1}\times r}$
and $C=[\mathbf{c}_{1}|\mathbf{c}_{2}|\ldots|\mathbf{c}_{r}]\in\mathbb{R}^{n_{2}\times r}$
and $\mathbf{x}\in\mathbb{R}^{n_{1}n_{2}}$ then 
\begin{eqnarray*}
\mathbf{y}=\left(B\odot C\right)^{T}\,\mathbf{x}=\left[\begin{array}{c}
\mathbf{c}_{1}^{T}X\mathbf{b}_{1}\\
\mathbf{c}_{2}^{T}X\mathbf{b}_{2}\\
\vdots\\
\mathbf{c}_{r}^{T}X\mathbf{b}_{r}
\end{array}\right].
\end{eqnarray*}
Where $X=\textrm{reshape}(\mathbf{x},n_{2},n_{1}).$\\
 \textbf{Proof:} Let $\mathbf{x}=[x_{1},x_{2},\ldots,x_{n_{1}n_{2}}]^{T}.$
Reshape the vector $\mathbf{x}$ as an $n_{2}\times n_{1}$ matrix
$X$ 
\begin{eqnarray*}
X=\left[\begin{array}{cccc}
x_{1} & x_{n_{2}+1} & ... & x_{(n_{1}-1)n_{2}+1}\\
x_{2} & x_{n_{2}+2} &  & x_{(n_{1}-1)n_{2}+2}\\
\vdots & \vdots\\
x_{n_{2}} & x_{2n_{2}} &  & x_{n_{1}n_{2}}
\end{array}\right].
\end{eqnarray*}
$(B\odot C)^{T}\mathbf{x}$ is given by 
\begin{eqnarray*}
(B\odot C)^{T}\mathbf{x}=[\mathbf{b}_{1}\otimes\mathbf{c}_{1}|\mathbf{b}_{2}\otimes\mathbf{c}_{2}|\ldots|\mathbf{b}_{r}
\otimes\mathbf{c}_{r}]_{n_{1}n_{2}\times r}^{T}\,\mathbf{x}_{n_{1}n_{2}\times1}\\
\\
=\left[\begin{array}{c}
(\mathbf{b}_{1}\otimes\mathbf{c}_{1})^{T}\\
(\mathbf{b}_{2}\otimes\mathbf{c}_{2})^{T}\\
\vdots\\
(\mathbf{b}_{r}\otimes\mathbf{c}_{r})^{T}
\end{array}\right]_{r\times n_{1}n_{2}}\mathbf{x}_{n_{1}n_{2}\times1}=\left[\begin{array}{c}
(\mathbf{b}_{1}\otimes\mathbf{c}_{1})^{T}\mathbf{x}\\
(\mathbf{b}_{2}\otimes\mathbf{c}_{2})^{T}\mathbf{x}\\
\vdots\\
(\mathbf{b}_{r}\otimes\mathbf{c}_{r})^{T}\mathbf{x}
\end{array}\right]\,\,=\,\left[\begin{array}{c}
\mathbf{c}_{1}^{T}X\mathbf{b}_{1}\\
\mathbf{c}_{2}^{T}X\mathbf{b}_{2}\\
\vdots\\
\mathbf{c}_{r}^{T}X\mathbf{b}_{r}
\end{array}\right].
\end{eqnarray*}
Where $X=\textrm{reshape}(\mathbf{x},n_{2},n_{1}).$ In the last step
of the above equation we have used 
$(\mathbf{b}_{i}^{T}\otimes\mathbf{c}_{i}^{T})\mathbf{x}=\mathbf{c}_{i}^{T}X\mathbf{b}_{i}.$
This can be shown easily in the following way.

Let $\mathbf{b}_{i}=[b_{1i},b_{2i},\ldots,b_{n_{1}i}]^{T}$ and 
$\mathbf{c}_{i}=[c_{1i},c_{2i},\ldots,c_{n_{2}i}]^{T}.$
Then 
\begin{eqnarray*}
(\mathbf{b}_{i}\otimes\mathbf{c}_{i})^{T}\mathbf{x}&=&\left[\begin{array}{c}
b_{1i}c_{1i}\\
b_{1i}c_{2i}\\
\vdots\\
b_{1i}c_{n_{2}i}\\
b_{2i}c_{1i}\\
\vdots\\
b_{n_{1}i}c_{n_{2}i}
\end{array}\right]^{T}\left[\begin{array}{c}
x_{1}\\
x_{2}\\
\vdots\\
x_{n_{2}}\\
x_{n_{2}+1}\\
\vdots\\
x_{n_{2}n_{1}}
\end{array}\right]\,\,\,\,\,\,\,\,\,\,\,\,\,\,\,\,\,\,\,\,\,\,\,\,\,\,\,\,\,\,\,\,\,\\
&=& b_{1i}c_{1i}x_{1}+b_{1i}c_{2i}x_{2}+........+b_{1i}c_{n_{2}i}x_{n_{2}}\\
&+& b_{2i}c_{1i}x_{n_{2}+1}+...................+b_{2i}c_{n_{2}i}x_{2n_{2}}\\
&\vdots& \\
&+& b_{n_{1}i}c_{1i}x_{n_{2}(n_{1}-1)+1}+........+b_{n_{1}i}c_{n_{2}i}x_{n_{1}n_{2}}
\end{eqnarray*}
\begin{eqnarray*}
&=& c_{1i}\left(b_{1i}x_{1}+b_{2i}x_{n_{2}+1}+.....+b_{n_{1}i}x_{n_{2}(n_{1}-1)+1}\right)\,\,\,\,\,\,\,\,\,\,\,\,\,\,\,\,\,
\,\,\,\,\,\,\,\,\,\,\,\,\,\,\,\,\,\,\,\,\,\,\\
&+& c_{2i}\left(b_{1i}x_{2}+b_{2i}x_{n_{2}+2}+.....+b_{n_{1}i}x_{n_{2}(n_{1}-1)+2}\right)\,\,\,\,\,\,\,\,\,\,\,\,\,\,\,\,\,
\,\,\,\,\,\,\,\,\,\,\,\,\,\,\,\,\,\,\,\,\,\,\\
&\vdots&\,\,\,\,\,\,\,\,\,\,\,\,\,\,\,\,\,\,\,\,\,\,\,\,\,\,\,\,\,\,\,\,\,\,\,\,\,\,\,\\
&+& c_{n_{2}i}\left(b_{1i}x_{n_{2}}+b_{2i}x_{2n_{2}}+..............+b_{n_{1}i}x_{n_{2}n_{1}}\right)\,\,\,\,\,\,\,\,\,\,\,\,
\,\,\,\,\,\,\,\,\,\,\,\,\,\,\,\,\,\,\,\,\,\,\,\,\,\,\,\\
\\
&=& [c_{1i},c_{2i},...,c_{n_{2}i}]\left[\begin{array}{c}
b_{1i}x_{1}+b_{2i}x_{n_{2}+1}+.....+b_{n_{1}i}x_{n_{2}(n_{1}-1)+1}\\
b_{1i}x_{2}+b_{2i}x_{n_{2}+2}+.....+b_{n_{1}i}x_{n_{2}(n_{1}-1)+2}\\
\vdots\\
b_{1i}x_{n_{2}}+b_{2i}x_{2n_{2}}+..............+b_{n_{1}i}x_{n_{2}n_{1}}
\end{array}\right]\\
&=&[c_{1i},c_{2i},...,c_{n_{2}i}]\left[\begin{array}{cccc}
x_{1} & x_{n_{2}+1} & .... & x_{(n_{1}-1)n_{2}+1}\\
x_{2} & x_{n_{2}+2} & .... & x_{(n_{1}-1)n_{2}+2}\\
\vdots & \vdots & \vdots\\
x_{n_{2}} & x_{2n_{2}} & ..... & x_{n_{1}n_{2}}
\end{array}\right]\left[\begin{array}{c}
b_{1i}\\
b_{2i}\\
\vdots\\
b_{n_{1}i}
\end{array}\right]\,\,\,\,\,\,\\
&=& \mathbf{c}_{i}^{T}X\mathbf{b}_{i}\,\,\, .
\end{eqnarray*}
Here $n_{2}(2n_{1}-1)$ operations are required to compute $X_{n_{2}n_{1}}\mathbf{b}_{i}$
and $(2n_{2}-1)$ operations to compute $\mathbf{c}_{i}^{T}(X_{n_{2}n_{1}}\mathbf{b}_{i}).$
Therefore the overall computational complexity is $2n_{1}n_{2}r+n_{2}r-r$
which is less than the complexity for computing $(B\odot C)^{T}\mathbf{x}$
directly.

\subsubsection{Generalization of Lemma 2}

Let us look at Lemma 2 in the case of three matrices 
$A_{4}\in\mathbb{R}^{n_{4}\times r},\,A_{3}\in\mathbb{R}^{n_{3}\times r}$
and $A_{2}\in\mathbb{R}^{n_{2}\times r}.$ Let $\mathbf{x}\in\mathbb{R}^{n_{4}n_{3}n_{2}}.$
We look at $(A_{4}\odot A_{3}\odot A_{2})^{T}\mathbf{x}.$ Reshape
the vector $\mathbf{x}$ into a matrix of size $n_{2}n_{3}\times n_{4}.$
Let $X_{4}=reshape(\mathbf{x},n_{2}n_{3},n_{4}).$ Then 
\begin{eqnarray*}
(A_{4}\odot A_{3}\odot A_{2})^{T}\mathbf{x}=\left[\begin{array}{c}
(A_{3}\odot A_{2})_{1}^{T}X_{4}(A_{4})_{1}\\
(A_{3}\odot A_{2})_{2}^{T}X_{4}(A_{4})_{2}\\
\vdots\\
(A_{3}\odot A_{2})_{r}^{T}X_{4}(A_{4})_{r}
\end{array}\right].
\end{eqnarray*}
Here $(A_{4})_{i}$ is the $i^{th}$ column of $A_{4}$ with size 
$n_{4}\times1.$ So the size of $X_{4}(A_{4})_{i}$
is $n_{2}n_{3}\times1.$ Let us denote the vector $X_{4}(A_{4})_{i}$
by 
\begin{eqnarray*}
(\mathbf{x}_{4})_{i}=X_{4}(A_{4})_{i},\,\,i=1,2,\ldots,r.
\end{eqnarray*}

Now reshape each $(\mathbf{x}_{4})_{i}$, $i=1,2,\ldots,r$ into matrices
$(X_{3})_{i}\in\mathbb{R}^{n_2\times n_3}.$ Then

\begin{eqnarray*}
(A_{4}\odot A_{3}\odot A_{2})^{T}\mathbf{x}=\left[\begin{array}{c}
(A_{3}\odot A_{2})_{1}^{T}(\mathbf{x}_{4})_{1}\\
(A_{3}\odot A_{2})_{2}^{T}(\mathbf{x}_{4})_{2}\\
\vdots\\
(A_{3}\odot A_{2})_{r}^{T}(\mathbf{x}_{4})_{r}
\end{array}\right]=\left[\begin{array}{c}
(A_{2})_{1}^{T}(X_{3})_{1}(A_{3})_{1}\\
(A_{2})_{2}^{T}(X_{3})_{2}(A_{3})_{2}\\
\vdots\\
(A_{2})_{r}^{T}(X_{3})_{r}(A_{3})_{r}
\end{array}\right].
\end{eqnarray*}

\subsubsection*{Computational complexity}

The computational complexity of the general product 
$(A_{L}\odot A_{L-1}\odot\ldots\odot A_{2})^{T}\mathbf{x}$ by 
the above technique is $r(2n-1)\left(\frac{n^{L-1}-1}{n-1}\right)$, 
where $A_{i}\in\mathbb{R}^{n\times r}$
and ${\mathbf{x}\in\mathbb{R}^{n^{L-1}}}$, whereas the direct computation
of this product is a bit more expensive, it requires $n^{L-1}r+O(rn^{L-1})$ arithmetic operations
including the computation of $A_{L}\odot A_{L-1}\odot\ldots\odot A_{2}$. 

\section{QCP Algorithm}

Let $f$ be a function discretized on a fine grid of 
size $2^{L}\,(\textrm{for example}\,L=15)$ with uniform length in an interval. 
The function values at the grid points generate a vector of size
$2^{L}.$ As described in the introduction we can reshape this long
vector as a tensor of order $L$ and one can approximate it as a sum
of products of vectors of length $2$. Fig. 1 shows an example of a
($3^{rd}$order tensor) quantized vector of length 
$2^{3},$ $[\tau_{1},\tau_{2},\tau_{3},\ldots,\tau_{8}]^{T}.$
The construction of a rank $r$ canonical approximation of such a 
tensor using Alternative Least Squares method is described below.

\begin{figure}[H]
\centering
~~~~~\includegraphics[scale=0.85]{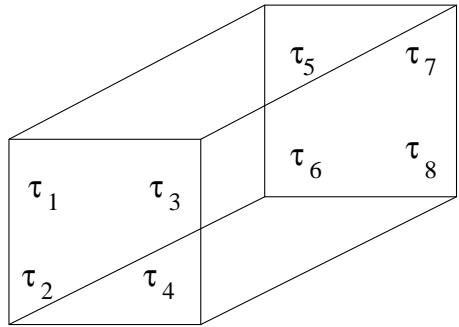}
\caption{$3^{rd}$ order tensor.}
\end{figure}

Let $I=[a,b].$ Consider an uniform mesh with mesh size $h=\frac{1}{2^{L}-1}.$
Let \textbf{$\mathbf{f}$} be the vector of length $2^{L}$ whose
entries are the values of the given function $f$ at these $2^{L}$
points on the grid. Let us denote $\mathbf{f}$ by 
\begin{eqnarray}
\mathbf{f}=[\tau_{1},\tau_{2},\tau_{3},\ldots,\tau_{2^{L}}]^{T}.\label{1}
\end{eqnarray}

Let $\chi$ be the quantized $L^{th}$ order tensor, given
by
\[\chi=reshape(\mathbf{f},\underset{L}{\underbrace{2,2,\ldots,2}})
\in\mathbb{R}^{2\times2\times...\times2}.\]
The precise definition of this operation is shortly recalled here: 

The vector $\mathbf{f}$ is reshaped to its quantics image in 
${\displaystyle {\displaystyle \otimes_{J=1}^{L}}}\mathbb{R}^{2}$ by dyadic folding, 
\begin{eqnarray*}
\mathcal{G}_{2,L}:\,\mathbf{f}\rightarrow\chi=\chi(\mathbf{j})\in\otimes_{J=1}^{L}\mathbb{R}^{2},\,\mathbf{j}=\{j_{1},j_{2},
\ldots,j_{L}\},\quad\text{with}\quad j_{\upsilon}\in\{1,2\},\upsilon=1,2,\ldots,L,
\end{eqnarray*}
 where for fixed $i,$ we have $\chi(\mathbf{j}):=\mathbf{f}(i)$
and $j_{\upsilon}=j_{\upsilon}(i)$ is defined via $2-$coding, $j_{\upsilon}-1=C_{-1+\upsilon},$
such that the coefficients $C_{-1+\upsilon}$ are found from the dyadic
representation of $i-1,$
\begin{eqnarray*}
i-1=C_{0}+C_{1}\,2+C_{2}\,2^{2}+\ldots+C_{L-1}\,2^{L-1}
\equiv{\displaystyle \sum_{\upsilon=1}^{L}(j_{\upsilon}-1)\,2^{\upsilon-1}.}
\end{eqnarray*}

The rank $r$ canonical approximation of the $L^{th}$ order tensor is 
\begin{eqnarray}
\chi\cong\sum_{k=1}^{r}\mathbf{a}_{k}^{(1)}\otimes\mathbf{a}_{k}^{(2)}\otimes\ldots\otimes\mathbf{a}_{k}^{(L)}\label{2}
\end{eqnarray}
where each $\mathbf{a}_{k}^{(i)}=\left[\begin{array}{c}
a_{1,k}^{(i)}\\
a_{2,k}^{(i)}
\end{array}\right]$ is a $2\times1$ vector and $\otimes$ is the usual tensor product.

Let 
\begin{eqnarray*}
A_{1}=[\mathbf{a}_{1}^{(1)},\mathbf{a}_{2}^{(1)},\ldots,\mathbf{a}_{r}^{(1)}],A_{2}=[\mathbf{a}_{1}^{(2)},
\mathbf{a}_{2}^{(2)},\ldots,\mathbf{a}_{r}^{(2)}],\ldots,A_{L}=[\mathbf{a}_{1}^{(L)},\mathbf{a}_{2}^{(L)},\ldots,\mathbf{a}_{r}^{(L)}].
\end{eqnarray*}
Here $A_{1},A_{2},\ldots,A_{L}$ are $2\times r$ matrices, corresponding
to $L$ different directions, whose columns are the unknown vectors
in equation (2).

The formulation of the ALS is the following: 
\begin{eqnarray}
\textrm{Minimize}\,\,\,\frac{1}{2}\left\Vert \chi-\sum_{k=1}^{r}\mathbf{a}_{k}^{(1)}\otimes\mathbf{a}_{k}^{(2)}\otimes
\ldots\otimes\mathbf{a}_{k}^{(L)}\right\Vert _{F}^{2},\label{3}
\end{eqnarray}
where $\left\Vert .\right\Vert _{F}$ is the Frobenius norm of a tensor.

In the ALS approach, this functional is minimized in an alternating
way. ALS fixes all $A_{j},\,j\neq i,j=1,2,\ldots,L$ to minimize for
$A_{i}$ and continue this process until some convergence criterion
is satisfied. That is, first fix $A_{2},A_{3},\ldots,A_{L}$ to solve
for $A_{1}$ and then fix $A_{1},A_{3},\ldots,A_{L}$ to solve for $A_{2}$
and so on and then fixes $A_{1},A_{2},\ldots,A_{L-1}$ to solve for $A_{L}$
and continue the process.

At each iteration of the ALS approach, we have $L$ steps. First,
we start with an initial guess on $A_{2},A_{3},\ldots,A_{L}$ and solve
for $A_{1}$, this gives the initial guess for the next step. Since
we are fixing $L-1$ matrices and solving for one of the matrices
$A_{i},\,i=1,2,\ldots,L$ at each step of an iteration, the problem is reduced
to a linear least-squares problem. \\
 In the $i^{th}$ step of an iteration, we fix $A_{1},A_{2},\ldots,A_{i-1},A_{i+1},\ldots,A_{L}$
and solve for\\
$A_{i}=\left[\begin{array}{cccc}
a_{1,1}^{(i)} & a_{1,2}^{(i)} & ... & a_{1,r}^{(i)}\\
a_{2,1}^{(i)} & a_{2,2}^{(i)} & ... & a_{2,r}^{(i)}
\end{array}\right].$ 
The resulting linear least-squares problem is: 
\begin{eqnarray}
\textrm{\textrm{minimize}}\,\,\mathcal{F},\,\,\,\,\,\,\,\,\,\,\,\,\,\,\,\,\,\,\,\,\,\,\,\,\,\,\,\,\,\,\,\,\,\,\,\,\,\,\,\,\,\,\,\,\,
\,\,\,\,\,\,\,\,\,\,\,\,\,\,\,\,\,\,\,\,\,\,\,\,\,\,\,\,\,\,\,\,\,\,\,\,\,\,\,\,\,\,\,\,\,\,\,\,\,\,\,\,\,\,\,\,\nonumber \\
\textrm{where}\,\,\mathcal{F}=\frac{1}{2}\left\Vert \chi-\sum_{k=1}^{r}\mathbf{a}_{k}^{(1)}\otimes\mathbf{a}_{k}^{(2)}
\otimes\ldots\otimes\mathbf{a}_{k}^{(L)}\right\Vert _{F}^{2}\,\,\textrm{with}\,\,A_{1},A_2,\ldots,A_{i-1},A_{i+1},
\ldots,A_{L}\,\,\textrm{fixed}.\label{4}
\end{eqnarray}

This gives the equations 
\begin{eqnarray*}
\frac{\partial\mathcal{F}}{\partial a_{1,1}^{(i)}}=\frac{\partial\mathcal{F}}{\partial a_{1,2}^{(i)}}=\ldots=
\frac{\partial\mathcal{F}}{\partial a_{1,r}^{(i)}}=0\quad\textrm{and}\quad\frac{\partial\mathcal{F}}{\partial a_{2,1}^{(i)}}=
\frac{\partial\mathcal{F}}{\partial a_{2,2}^{(i)}}=\ldots=\frac{\partial\mathcal{F}}{\partial a_{2,r}^{(i)}}=0.
\end{eqnarray*}

These equations can be written in the form 
\begin{eqnarray}
\left[\begin{array}{cc}
\hat{A}_{i}^{T}\hat{A}_{i}\\
 & \hat{A}_{i}^{T}\hat{A}_{i}
\end{array}\right]\left[\begin{array}{c}
\left[\begin{array}{c}
a_{1,1}^{(i)}\\
a_{1,2}^{(i)}\\
\vdots\\
a_{1,r}^{(i)}
\end{array}\right]\\
\begin{array}{c}
\left[\begin{array}{c}
a_{2,1}^{(i)}\\
a_{2,2}^{(i)}\\
\vdots\\
a_{2,r}^{(i)}
\end{array}\right]\end{array}
\end{array}\right]=\left[\begin{array}{c}
\hat{A}_{i}^{T}\chi(:,\ldots,:,1,:,\ldots,:)\\
\hat{A}_{i}^{T}\chi(:,\ldots,:,2,:,\ldots,:)
\end{array}\right].\label{5}
\end{eqnarray}
Here $\hat{A}_{i}=A_{L}\odot A_{L-1}\odot...\odot A_{i+1}\odot A_{i-1}\odot...\odot A_{1},$
where $\odot$ denotes the Khatri-Rao product of matrices (see A2). This
$\hat{A}_{i}$ is a $2^{L-1}\times r$ matrix and $\chi(:,\ldots,:,j,:,\ldots,:)$~~$j=1,2$
is a vector of length $2^{L-1}.$ $\hat{A}_{i}^{T}\hat{A}_{i}$ is a
$r\times r$ symmetric positive definite matrix. The direct computation
of $\hat{A}_{i}^{T}\hat{A}_{i}$ and $\hat{A}_{i}^{T}\chi(:,\ldots,:,j,:,\ldots,:)$
is expensive. The fast computation of these products are described in section 2.\\
\medskip{}
\textbf{Remark:} For a better understanding of the structure of $\hat{A}_{i}$
and the derivation of (5) we refer to appendix A3. All steps
of the ALS algorithm for rank-2 canonical approximation of a $4^{th}$
order tensor are shown there in detail. 

>From equation (5) one can see that we need to solve two $r\times r$
linear systems with the same matrix $\hat{A}_{i}^{T}\hat{A}_{i}$
and different right hand side vectors at each step of an iteration.
We continue the iterations until a convergence criterion is reached.

\subsubsection*{Algorithm}

\textemdash \textemdash \textemdash \textemdash \textemdash \textemdash \textemdash \textemdash \textemdash 
\textemdash \textemdash \textemdash \textemdash \textemdash \textemdash \textemdash \textemdash \textemdash 
\textemdash \textemdash \textemdash \textemdash \textemdash \textemdash \textemdash \textemdash \textemdash 
\textemdash \textemdash \textemdash \textemdash \textemdash \textemdash \textemdash \textemdash \textemdash 
\textemdash \textemdash \textemdash \textemdash \textemdash{}

\textit{Define tolerance $\epsilon$}

\textit{Maximum iterations=Maxiter}

\textit{Initialize $A_{i}\in\mathbb{R}^{2\times r},i=1,2,\ldots,L.$ }

\textit{while iter<=Maxiter}

\textit{$C_{i}=A_{i},i=1,2,\ldots,L$}

\textit{\medskip{}
 }

\textit{for $i=1,2,..L$}

\textit{Obtain $\hat{A}_{i}=A_{L}\odot A_{L-1}\odot\ldots\odot A_{i+1}\odot A_{i-1}\odot\ldots\odot A_{1}$;
$\hat{A}_{i}^{T}\chi(:,\ldots,:,j,:,\ldots,:)$ for $j=1,2$}

\textit{Solve $\hat{A}_{i}^{T}\hat{A}_{i}$$\left[\begin{array}{c}
a_{1,1}^{(i)}\\
a_{1,2}^{(i)}\\
\vdots\\
a_{1,r}^{(i)}
\end{array}\right]=\hat{A}_{i}^{T}\chi(:,\ldots,:,1,:,\ldots,:)\,\,\textrm{and analog}\,\,\,\hat{A}_{i}^{T}\hat{A}_{i}\begin{array}{c}
\left[\begin{array}{c}
a_{2,1}^{(i)}\\
a_{2,2}^{(i)}\\
\vdots\\
a_{2,r}^{(i)}
\end{array}\right]\end{array}$}

\textit{end for}

\textit{\medskip{}
 }

\textit{stop if $max\{max\left|A_{i}-C_{i}\right|\}<\epsilon$ }

\textit{iter=iter+1}

\textit{end while}

\textemdash \textemdash \textemdash \textemdash \textemdash \textemdash \textemdash \textemdash \textemdash 
\textemdash \textemdash \textemdash \textemdash \textemdash \textemdash \textemdash \textemdash \textemdash 
\textemdash \textemdash \textemdash \textemdash \textemdash \textemdash \textemdash \textemdash \textemdash 
\textemdash \textemdash \textemdash \textemdash \textemdash \textemdash \textemdash \textemdash \textemdash 
\textemdash \textemdash \textemdash \textemdash \textendash{}

\subsubsection*{Computational Complexity}

Let the number of iterations in the above algorithm be $iter.$ In
each iteration step of ALS we need to compute 
$\hat{A}_{i}=A_{L}\odot A_{L-1}\odot\ldots\odot A_{i+1}\odot A_{i-1}\odot\ldots\odot A_{1}$
and $\hat{A}_{i}^{T}\chi(:,\ldots,:,1,:,\ldots,:)$ as well as 
$\hat{A}_{i}^{T}\chi(:,\ldots,:,2,:,\ldots,:)$
for $i=1,2,\ldots,L$ and need to solve a linear least-squares system twice. 

The computation $\hat{A}_{i}=A_{L}\odot A_{L-1}\odot\ldots\odot A_{i+1}\odot A_{i-1}\odot\ldots\odot A_{1}$
requires $O((L-1)2r^{2})$ arithmetic operations (look at section
2, and here $n=2$) and $\hat{A}_{i}^{T}\chi(:,\ldots,:,j,:,\ldots,:)$ requires
$3r(2^{L-1}-1)$ operations (look at section 2 ). $O(r^{3})$ operations
are required to solve a $r\times r$ linear system. So, the total
complexity of the algorithm is $O\left(\,L\,\left((L-1)2r^{2}+3r(2^{L-1}-1)+r^{3}\right)\right)$
for each iteration step. That is $O(2^{L-1})$ per iteration
step. 

\subsubsection*{Comments on the algorithm}

This is a straight forward ALS algorithm applied to higher order
tensors of order $L$. The initialization of 
\textit{$A_{i}\in\mathbb{R}^{2\times r},i=1,2,\ldots,L$}{\small{}
is random}\textit{.} The condition number of the matrices $\hat{A}_{i}^{T}\hat{A}_{i}$
is large for large values of $r.$ 

\section{Numerical examples}

In this section we present the canonical approximation of some functions
discretized on $[0,1]$ and consider the approximation in the following
format 
\begin{eqnarray*}
\sum_{k=1}^{r}\left(\begin{array}{c}
1\\
a_{2,k}^{(1)}
\end{array}\right)\otimes\left(\begin{array}{c}
1\\
a_{2,k}^{(2)}
\end{array}\right)\otimes\ldots\otimes\left(\begin{array}{c}
a_{1,k}^{(L)}\\
a_{2,k}^{(L)}
\end{array}\right).
\end{eqnarray*}

The number of parameters in this format is almost half of the parameters
required for the canonical representation given in equation (2). So,
the computational complexity is here further reduced. The condition
numbers of the matrices $\hat{A}_{i}^{T}\hat{A}_{i}$ are much better
in this case. 

In all the numerical examples given below, the functions are discretized
on a uniform grid of size $2^{15},$ so the reshaped tensor is of
order $15.$ In all the tables below ``$error"$ denotes the maximum error in the canonical
approximation of the discretized function. The initial matrices $A_{i}$
are chosen randomly and the computations are carried out in MATLAB.
\medskip{}

\begin{table}[H]
~~~~~~~~~~~~~~~~~~~~~~~~~~~~~~~~~~~~~~~~~~~~~~~~~{\small{}~~~~~~~~}%
\begin{tabular}{|c|c|}
\hline 
{\small{}$r$ } & {\small{}$error$}\tabularnewline
\hline 
\hline 
{\small{}1 } & {\small{}0.108596}\tabularnewline
\hline 
{\small{}2 } & {\small{}0.031}\tabularnewline
\hline 
{\small{}3 } & {\small{}0.0081}\tabularnewline
\hline 
{\small{}4 } & {\small{}0.0023}\tabularnewline
\hline 
{\small{}5 } & {\small{}0.00071}\tabularnewline
\hline 
{\small{}6 } & {\small{}0.00024}\tabularnewline
\hline 
{\small{}7 } & {\small{}0.00015}\tabularnewline
\hline 
{\small{}8 } & {\small{}0.0000881}\tabularnewline
\hline 
{\small{}9 } & {\small{}0.0000461}\tabularnewline
\hline 
{\small{}10 } & {\small{}0.0000210}\tabularnewline
\hline 
\end{tabular}{\small \par}

\protect\caption{Error in the maximum norm for different values of $r$. }
\end{table}

\textbf{Example 1:} Consider the function $f(x)=e^{-x^{2}}$ in $[0,1]$.
We have obtained the canonical approximation with different ranks, see Table 1.\\

\textbf{$\!\!\!\!\!\!\!\!\!$Example 2:}
Consider the functions $sin(\pi x),sin(2\pi x),\,\,\textrm{and}\,\,sin(4\pi x)$
in $[0,1]$. Table 2 shows the error in the maximum norm for different
values of $r.$

\begin{table}[H]
~~~~~~~~~~~~~~~~~~~~~~~~~~~~~~~~~~~{\small{}~~~~~~}%
\begin{tabular}{|c|c|cc|cc|}
\hline 
\multicolumn{1}{|c}{} & \multicolumn{1}{c}{{\small{}$sin(\pi x)$}} & \multicolumn{1}{c}{} & 
\multicolumn{1}{c}{{\small{}$sin(2\pi x)$}} & \multicolumn{1}{c}{} & {\small{}$sin(4\pi x)$}\tabularnewline
\hline 
{\small{}$r$ } & {\small{}$error$ } &  & {\small{}$errror$ } &  & {\small{}$error$}\tabularnewline
\hline 
\hline 
{\small{}1 } & {\small{}0.63658 } &  & {\small{}1.000 } &  & {\small{}1.0}\tabularnewline
\hline 
{\small{}2 } & {\small{}0.164 } &  & {\small{}0.250 } &  & {\small{}0.162}\tabularnewline
\hline 
{\small{}3 } & {\small{}0.0336 } &  & {\small{}0.0723 } &  & {\small{}0.067}\tabularnewline
\hline 
{\small{}4 } & {\small{}0.00635 } &  & {\small{}0.0341 } &  & {\small{}0.0308}\tabularnewline
\hline 
{\small{}5 } & {\small{}0.0014 } &  & {\small{}0.00591 } &  & {\small{}0.0059}\tabularnewline
\hline 
{\small{}6 } & {\small{}0.000292 } &  & {\small{}0.00168 } &  & {\small{}0.0022}\tabularnewline
\hline 
{\small{}7 } & {\small{}0.0000822 } &  & {\small{}0.000389 } &  & {\small{}0.0010}\tabularnewline
\hline 
{\small{}8 } & {\small{}0.0000572 } &  & {\small{}0.000172 } &  & {\small{}0.000370}\tabularnewline
\hline 
{\small{}9 } & {\small{}0.00000901 } &  & {\small{}0.0000886 } &  & {\small{}0.000142}\tabularnewline
\hline 
{\small{}10 } & {\small{}0.00000671 } &  & {\small{}0.0000317 } &  & {\small{}0.000070}\tabularnewline
\hline 
\end{tabular}{\small \par}

\protect\caption{Error for different ranks in the canonical approximation.}
\end{table}
\textbf{$\!\!\!\!\!\!\!\!\!$Example 3:} Now consider the functions
$f(x)=x$ or $f(x)=x^{2}$ in $[0,1]$. Table 3 shows the $error$ for different
values of $r.$\\

\begin{table}[H]
~~~~~~~~~~~~~~~~~~~~~~~~~~~~~~~~~~~~~~~{\small{}~~~~~}%
\begin{tabular}{|c|c|c|}
\hline 
\multicolumn{1}{|c}{} & \multicolumn{1}{c}{{\small{}$x$}} & {\small{}$x^{2}$}\tabularnewline
\hline 
{\small{}$r$ } & {\small{}$error$ } & {\small{}$error$}\tabularnewline
\hline 
\hline 
{\small{}1 } & {\small{}0.176 } & {\small{}0.075}\tabularnewline
\hline 
{\small{}2 } & {\small{}0.0186 } & {\small{}0.0276}\tabularnewline
\hline 
{\small{}3 } & {\small{}0.00576 } & {\small{}0.00661}\tabularnewline
\hline 
{\small{}4 } & {\small{}0.00133 } & {\small{}0.00121}\tabularnewline
\hline 
{\small{}5 } & {\small{}0.000346 } & {\small{}0.000218}\tabularnewline
\hline 
{\small{}6 } & {\small{}0.000082 } & {\small{}0.00005}\tabularnewline
\hline 
{\small{}7 } & {\small{}0.000022 } & {\small{}0.0000125}\tabularnewline
\hline 
{\small{}8 } & {\small{}0.00000652 } & {\small{}0.00000927}\tabularnewline
\hline 
{\small{}9 } & {\small{}0.00000268 } & {\small{}0.00000351}\tabularnewline
\hline 
{\small{}10 } & {\small{}0.000000728 } & {\small{}0.00000252}\tabularnewline
\hline 
\end{tabular}{\small \par}

\protect\caption{Error for different values of $r$.}
\end{table}

In all the examples above, one can observe that the error decays exponentially
with $r$, like $\mu^{r}$ where $\mu<1.$ Also,
one can see that the function (or better: its discretized representation) 
has been well approximated by 
the QCP format using only $160$ parameters, where the original size was $2^{15}$. 
Please note that so 
far we have used complete information of the data to obtain the QCP approximation. 
A more effective way based on the QCP interpolation is sketched in the following section.

\section{The QCP approximation using only a few function calls}

In section 3, we have seen the construction of a rank $r$ canonical
approximation using the complete data of size $2^{L}.$ Here we describe
the idea of constructing the rank $r$ canonical
approximation using function values at a few sampling points only.
The more detailed presentation is the topic of our ongoing work.
Let $M\,(=O(2Lr))$ be the number of sampling points, comparable to the number of unknown representation 
parameters. 
Many issues like a good choice of the sampling points and the robust error analysis of the method 
will not be addressed in this article. This approach can be viewed as the sparse interpolation
of a given function in the QCP format by using a small number of functional calls. 

Consider the rank-$r$ canonical approximation of the tensor $\chi$
\begin{eqnarray*}
\chi\cong\sum_{k=1}^{r}\mathbf{a}_{k}^{(1)}\otimes
\mathbf{a}_{k}^{(2)}\otimes\ldots\otimes\mathbf{a}_{k}^{(L)}.
\end{eqnarray*}
The method to evaluate the unknown parameters $\mathbf{a}_{k}^{(i)},i=1,2,\ldots,L,\;k=1,2,\ldots,r$,
using the information of the tensor $\chi$ only at $M$ positions
is given below.
 We let 
\begin{eqnarray*}
A_{1}=[\mathbf{a}_{1}^{(1)},\mathbf{a}_{2}^{(1)},\ldots,\mathbf{a}_{r}^{(1)}],A_{2}
=[\mathbf{a}_{1}^{(2)},
\mathbf{a}_{2}^{(2)},\ldots,\mathbf{a}_{r}^{(2)}],\ldots,A_{L}
=[\mathbf{a}_{1}^{(L)},\mathbf{a}_{2}^{(L)},\ldots,\mathbf{a}_{r}^{(L)}]
\end{eqnarray*}
be the side matrices.

Suppose we haven chosen $M$ points $s_{1},s_{2},\ldots,s_{M}$
on the grid with corresponding function values such
that they represent the function well in the whole interval. 
The corresponding entries in the vector $\mathbf{f}$ are 
denoted by $\tau_{s_{1}},\tau_{s_{2}},\ldots,\tau_{s_{M}}.$
We can identify these entries at certain positions in the $L^{th}$
order tensor $\chi$ and one can obtain the subscripts in the tensor
product corresponding to the linear index of the entries 
$\tau_{s_{1}},\tau_{s_{2}},\ldots,\tau_{s_{M}}.$
Let us denote the subscripts corresponding to each linear index by
\begin{eqnarray}
s_{1}&\rightarrow&\left(i_{1}^{s_{1}},i_{2}^{s_{1}},\ldots,i_{L}^{s_{1}}\right)\nonumber \\
s_{2}&\rightarrow&\left(i_{1}^{s_{2}},i_{2}^{s_{2}},\ldots,i_{L}^{s_{2}}\right)\nonumber \\
&\vdots&\,\,\,\,\,\,\,\,\,\label{eq:2}\\
s_{M}&\rightarrow&\left(i_{1}^{s_{M}},i_{2}^{s_{M}},\ldots,i_{L}^{s_{M}}\right).\nonumber 
\end{eqnarray}

Remember that each subscript $i_{j}^{s_{k}}$ is either $1$ or $2$ for all
$j=1,2,\ldots,L,\;k=1,2,\ldots,M.$ 

Analog to what is shown in Appendix A3, we minimize the functional 
\begin{eqnarray*}
\mathcal{F}=\frac{1}{2}\left(\left(\tau_{s_{1}}-{\displaystyle \sum_{k=1}^{r}}a_{i_{1}^{s_{1}},k}^{(1)}\,
\,a_{i_{2}^{s_{1}},k}^{(2)}\cdots a_{i_{L}^{s_{1}},k}^{(L)}\right)^{2}+\left(\tau_{s_{2}}-{\displaystyle 
\sum_{k=1}^{r}}a_{i_{1}^{s_{2}},k}^{(1)}\,\,a_{i_{2}^{s_{2}},k}^{(2)}\cdots 
a_{i_{L}^{s_{2}},k}^{(L)}\right)^{2}\right.\\
+\ldots+\left.\left(\tau_{s_{M}}-{\displaystyle 
\sum_{k=1}^{r}}a_{i_{1}^{s_{M}},k}^{(1)}\,\,a_{i_{2}^{s_{M}},k}^{(2)}
\cdots a_{i_{L}^{s_{M}},k}^{(L)}\right)^{2}\right) \to \min
\end{eqnarray*}
with respect to the unknown side matrices.
At each iteration of ALS we have $L$ steps. In the $i^{th}$ step
of an iteration, we fix $A_{1},A_{2},\ldots,A_{i-1},$ $A_{i+1},\ldots,A_{L}$
and solve for{\small{} $A_{i}=\left[\begin{array}{cccc}
a_{1,1}^{(i)} & a_{1,2}^{(i)} & ... & a_{1,r}^{(i)}\\
a_{2,1}^{(i)} & a_{2,2}^{(i)} & ... & a_{2,r}^{(i)}
\end{array}\right].$} This reduces the problem to a linear least-squares problem. The
linear system looks very similar to the system in (7) but with some differences. 
Here we describe it in detail.

Among the $M$ sampling points $s_{1},s_{2},\ldots,s_{M},$ let $p_{1},p_{2},\ldots,p_{N_{1}}$ be the
linear indices having $1$ as the $i^{th}$ subscript and $q_{1},q_{2},\ldots,q_{N_{2}}$
be the linear indices having 2 as the $i^{th}$ subscript ($N_{1}+N_{2}=M)$.
Then the linear system is given by 
\begin{eqnarray*}
\hat{A}_{i,1}^{T}\hat{A}_{i,1}\left[\begin{array}{c}
a_{1,1}^{(i)}\\
a_{1,2}^{(i)}\\
\vdots\\
a_{1,r}^{(i)}
\end{array}\right]=\hat{A}_{i,1}^{T}r_{i,1}\,\,\textrm{and}\,\,
\hat{A}_{i,2}^{T}\hat{A}_{i,2}\begin{array}{c}
\left[\begin{array}{c}
a_{2,1}^{(i)}\\
a_{2,2}^{(i)}\\
\vdots\\
a_{2,r}^{(i)}
\end{array}\right]\end{array}=\hat{A}_{i,2}^{T}r_{i,2},
\end{eqnarray*}

\begin{eqnarray*}
\text{where}\quad\hat{A}_{i,1}&=&\left[\begin{array}{cccc}
\hat{a}_{1}^{p_{1}} & \hat{a}_{2}^{p_{1}} &  & \hat{a}_{r}^{p_{1}}\\
\hat{a}_{1}^{p_{2}} & \hat{a}_{2}^{p_{2}} &  & \hat{a}_{r}^{p_{2}}\\
\\
\hat{a}_{1}^{p_{N_{1}}} & \hat{a}_{2}^{p_{N_{1}}} &  & \hat{a}_{r}^{p_{N_{1}}}
\end{array}\right]\,\,\textrm{with}\,\,\,\hat{a}_{k}^{p}=a_{i_{L}^{p},k}^{(L)}\,a_{i_{L-1}^{p},k}^{(L-1)}
\cdots a_{i_{i+1}^{p},k}^{(i+1)}\,a_{i_{i-1}^{p},k}^{(i-1)}\cdots a_{i_{1}^{p},k}^{(1)}\,\,,\\
\hat{A}_{i,2}&=&\left[\begin{array}{cccc}
\hat{a}_{1}^{q_{1}} & \hat{a}_{2}^{q_{1}} &  & \hat{a}_{r}^{q_{1}}\\
\hat{a}_{1}^{q_{2}} & \hat{a}_{2}^{q_{2}} &  & \hat{a}_{r}^{q_{2}}\\
\\
\hat{a}_{1}^{q_{N_{2}}} & \hat{a}_{2}^{q_{N_{2}}} &  & \hat{a}_{r}^{q_{N_{2}}}
\end{array}\right]\,\,\textrm{with}\,\,\,\hat{a}_{k}^{q}=a_{i_{L}^{q},k}^{(L)}\,a_{i_{L-1}^{q},k}^{(L-1)}
\cdots a_{i_{i+1}^{q},k}^{(i+1)}\,a_{i_{i-1}^{q},k}^{(i-1)}\cdots a_{i_{1}^{q},k}^{(1)}\\
\text{and}\quad r_{i,1}&=&\left[\begin{array}{c}
\tau_{p_{1}}\\
\tau_{p_{2}}\\
\vdots\\
\tau_{p_{N_{1}}}
\end{array}\right]\,\,,r_{i,2}=\left[\begin{array}{c}
\tau_{q_{1}}\\
\tau_{q_{2}}\\
\vdots\\
\tau_{q_{N_{2}}}
\end{array}\right].
\end{eqnarray*}

\textbf{\textit{Remark:}}\textit{ The matrices $\hat{A}_{i,1}\,\textrm{or}\,\hat{A}_{i,2}$
are very similar to $A_{L}\odot A_{L-1}\odot\ldots\odot A_{i+1}\odot A_{i-1}\odot\ldots\odot A_{1}$
but with many rows missing. The sizes of the matrices $\hat{A}_{i,1},\,\textrm{and}\,\,\hat{A}_{i,2}$
are $N_{1}\times r$ and $N_{2}\times r$ respectively, which are very
small compared to $\hat{A}_{i}$ in (5).} \medskip{}

This leads to a reduction of the computational complexity. Here we present a numerical
example to check the performance of the algorithm. 
We consider an approximation in the following format
\begin{eqnarray*}
\sum_{k=1}^{r}\left(\begin{array}{c}
a_{1,k}^{(1)}\\
a_{2,k}^{(1)}
\end{array}\right)\otimes\left(\begin{array}{c}
a_{1,k}^{(2)}\\
a_{2,k}^{(2)}
\end{array}\right)\otimes....\otimes\left(\begin{array}{c}
a_{1,k}^{(L)}\\
a_{2,k}^{(L)}
\end{array}\right).
\end{eqnarray*}

A further reduction of the number of unknowns is possible if one uses 
the format which has been discussed in section 4. \medskip{}
\\
\textbf{Example 4:} Consider the function $f(x)=e^{-x^{2}}$ in $[0,1]$
and $f(x)=e^{-50x^{2}}$ in $[0,0.25].$ Let $L=12$ and therefore
the grid size is $2^{12}.$ We have obtained the canonical approximation
with different ranks using the information of the function at $M=2Lr$
or $M=4Lr$ sampling points. The sampling points and initial matrices
are chosen randomly. Table 4 shows ``$error$'' in the approximation
for different values of the rank $r$ (in analogy to section 4, the maximum error is considered). 

\begin{table}[H]
~~~~~~~~~~~~%
\begin{tabular}{|c|c|c|c|c|c|c|}
\hline 
\multicolumn{1}{|c}{} & \multicolumn{1}{c}{} & \multicolumn{1}{c}{$e^{-x^{2}}$ in $[0,1]$} & 
\multicolumn{1}{c}{} & 
& \multicolumn{1}{c}{~~~$e^{-50x^{2}}$ in~~$[0,0.25]$} & \tabularnewline
\hline 
$r$ & $M=2Lr$ & $error$ & $M=4Lr$ & $error$ & $M=4Lr$ & $error$\tabularnewline
\hline 
\hline 
1 & 24 & 0.219347 & 48 & 0.144140 & 48 & 0.2081219\tabularnewline
\hline 
2 & 48 & 0.056676 & 96 & 0.0291372 & 96 & 0.0291072\tabularnewline
\hline 
3 & 72 & 0.011712 & 144 & 0.0075389 & 144 & 0.0124090\tabularnewline
\hline 
4 & 96 & 0.006980 & 192 & 0.0036845 & 192 & 0.0040713\tabularnewline
\hline 
5 & 120 & 0.003715 & 240 & 0.0019918 & 240 & 0.0023895\tabularnewline
\hline 
6 & 144 & 0.002515 & 288 & 0.0002400 & 288 & 0.0013455\tabularnewline
\hline 
7 & 168 & 0.001142 &  &  & 336 & 0.00084574\tabularnewline
\hline 
8 & 192 & 0.000697 &  &  & 384 & 0.00026631\tabularnewline
\hline 
\end{tabular}\caption{Error of the QCP interpolation for different values of $r$ and $M$.}
\end{table}

Table 4 also shows the number of sampling points used to obtain the
canonical approximation. For the function $f(x)=e^{-x^{2}},$ the
$error$ decays very fast in the case of $M=4Lr$ compared to the
case of $M=2Lr.$ One can see that we have used function values
only at $288$ points to approximate the tensor to $O(10^{-4})$ accuracy
instead of using the information at $4096$\textbf{ }points. The results
are presented for $M=4Lr$ in the case of the sharp Gaussian $f(x)=e^{-50x^{2}}.$
The error decays fast and we have used the information only at $384$
points to approximate the tensor to $O(10^{-4})$ instead of $4096.$
In both cases, we can see that the error decays exponentially like $\mu^{r},$
where $\mu<1.$ 

\begin{figure}[H]
\centering
\includegraphics[scale=0.75]{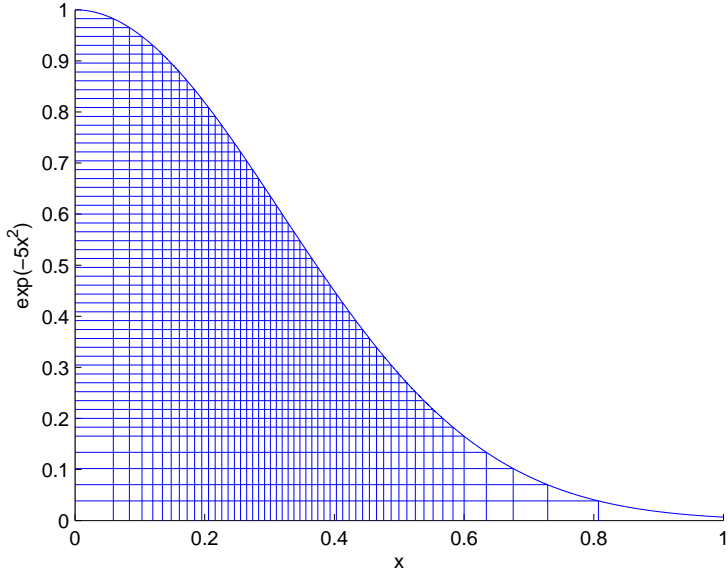}
\caption{Adaptive grid for the sparse QCP interpolation.}
 \label{Adapt_QTT_interp}
\end{figure}

The sparse interpolation in the QCP format requires the information of the function
only at $M\,(\sim 2Lr)$ points instead of the information at the full set $2^{L}$ of grid 
points. The overall computational complexity of the algorithm is reduced dramatically
and it is $O(M)$ per iteration step of the ALS algorithm. Here
$M \ll N=2^{L}.$ In the above numerical example 
the sampling points were chosen randomly. Clearly, there are many strategies for 
adaptive selection of sampling points 
based on some a priori knowledge about the behavior 
of the underlying function, but this issue will not be discussed here in detail. 
Figure \ref{Adapt_QTT_interp} shows an example of the adaptive choice of the interpolation grid
for the function $e^{-5 x^2}$.
Notice that the so-called TT-cross approximation in the TT format \cite{OsTy_TT2:09}
requires asymptotically smaller number of functional calls than $N$ in the case of large enough $N$.

\section{Conclusions and future work}

%

In this article, the ALS-type algorithms for 
approximation/interpolation of a function in QCP format have been described. 
The representation complexity of the rank-$r$ QCP format is estimated by $2Lr$.
As commented in section 3, the condition numbers of
the matrices appearing in each iteration of the ALS algorithm are large
for large values of the rank $r.$ Complete data of the tensor
has been used to obtain the CP approximation at the computational
cost $O(2^{L-1})$ per iteration. This complexity is reduced
if the approximation can be obtained using only a few data points, which
can be viewed as the sparse interpolation of a given function in the QCP
format. 

The idea of obtaining CP approximation using only small
number of functional calls is described and numerical examples are presented. 
In this case the overall computational complexity of the QCP approximation is
only $O(2Lr)$ per iteration step of the algorithm, i.e., it is proportional to the 
number of representation parameters in the target QCP tensor. It is remarkable that the complexity
of the QCP interpolation scales linearly in the CP rank and logarithmically in the 
full vector size.

A discussion of different strategies for clever choice of the sampling points as well as the error 
analysis of the method and the
extension of the algorithm to functions of two or three variables
is postponed to ongoing work. The QCP format can be used in the approximation of the
solution of PDEs, integration of highly oscillating functions and to
approximate functions where the calculation of function values
is computationally expensive. This format can also be used to just represent
functions that depend on many parameters. 

\vspace{0.3cm}
{\bf Acknowledgements.} 
KKN appreciates the support provided by the Max-Planck Institute
for Mathematics in the Sciences (Leipzig, Germany) during his scientific visit in 2015.
The authors are thankful to Dr. V. Khoromskaia (MPI MIS, Leipzig) for useful discussions.

\section*{Appendix}

\subsection*{A1. Kronecker product}

Let $A$ be an $m\times n$ matrix and $B$ be an $p\times q$ matrix, then
the Kronecker product $A\otimes B$ is an $mp\times nq$ block matrix:
\begin{eqnarray*}
A\otimes B=\left[\begin{array}{cccc}
a_{11}B & a_{12}B & ... & a_{1n}B\\
a_{21}B & a_{22}B & ... & a_{2n}B\\
\\
a_{m1}B & a_{m2}B & ... & a_{mn}B
\end{array}\right].
\end{eqnarray*}

\subsubsection*{Properties of the Kronecker product}

For matrices $B,C,D$ and $F$ of suitable sizes the following properties hold:

P1. $(B\otimes C)^{T}=B^{T}\otimes C^{T}.$\\
 P2. $(B\otimes C)\otimes D=B\otimes(C\otimes D).$\\
 P3. $(B\otimes C)(D\otimes F)=BD\otimes CF.$\\
 P4. $(B\otimes C)^{-1}=B^{-1}\otimes C^{-1}.$

\subsection*{A2. Khatri-Rao product}

\subsubsection*{\textmd{Let $B=[\mathbf{b}_{1}|\mathbf{b}_{2}|\ldots|\mathbf{b}_{r}]\in\mathbb{R}^{n_{1}\times r},$
where $\mathbf{b}_{1},\mathbf{b}_{2},\ldots,\mathbf{b}_{r}$ are the columns
of the matrix $B.$ Let $C=[\mathbf{c}_{1}|\mathbf{c}_{2}|\ldots|\mathbf{c}_{r}]\in\mathbb{R}^{n_{2}\times r}.$
The Khatri-Rao product of $B$ and $C$ is defined as the $n_{1}n_{2}\times r$ matrix
\begin{eqnarray*}
B\odot C=[\mathbf{b}_{1}\otimes\mathbf{c}_{1}|\mathbf{b}_{2}\otimes\mathbf{c}_{2}|\ldots|\mathbf{b}_{r}\otimes\mathbf{c}_{r}].
\end{eqnarray*}
Here $\mathbf{b}_{i}\times\mathbf{c}_{i}$
is} 
\begin{eqnarray*}
\mathbf{b}_{i}\otimes\mathbf{c}_{i}=\left[\protect\begin{array}{c}
b_{1i}\protect\\
b_{2i}\protect\\
\vdots\protect\\
b_{n_{1}i}
\protect\end{array}\right]\otimes c_{i}=\left[\protect\begin{array}{c}
b_{1i}c_{i}\protect\\
b_{2i}c_{i}\protect\\
\vdots\protect\\
b_{n_{1}i}c_{i}
\protect\end{array}\right]=\left[\protect\begin{array}{c}
b_{1i}c_{1i}\protect\\
b_{1i}c_{2i}\protect\\
\vdots\protect\\
b_{1i}c_{n_{2}i}\protect\\
b_{2i}c_{1i}\protect\\
\vdots\protect\\
b_{n_{1}i}c_{n_{2}i}
\protect\end{array}\right].
\end{eqnarray*}
}

$O(n_{1}n_{2}r)$ arithmetic operations are required to compute $B\odot C.$

\subsection*{A3. Rank 2 canonical approximation of a $4^{th}$ order tensor}

Consider a $4^{th}$ order tensor. Imagine that the tensor is generated
by reshaping a vector $\mathbf{x}=[\tau_{1},\,\tau_{2},\,\tau_{3},\ldots,\tau_{16}]^{T}$
of length 16. Let $\chi=reshape(\mathbf{x},2,2,2,2).$

We obtain a rank two canonical approximation to the tensor $\chi$
using ALS. The rank-2 approximation in canonical format is given by
\begin{eqnarray*}
\chi\cong\left(\begin{array}{c}
a_{1}^{1}\\
a_{2}^{1}
\end{array}\right)\otimes\left(\begin{array}{c}
b_{1}^{1}\\
b_{2}^{1}
\end{array}\right)\otimes\left(\begin{array}{c}
c_{1}^{1}\\
c_{2}^{1}
\end{array}\right)\otimes\left(\begin{array}{c}
d_{1}^{1}\\
d_{2}^{1}
\end{array}\right)+\left(\begin{array}{c}
a_{1}^{2}\\
a_{2}^{2}
\end{array}\right)\otimes\left(\begin{array}{c}
b_{1}^{2}\\
b_{2}^{2}
\end{array}\right)\otimes\left(\begin{array}{c}
c_{1}^{2}\\
c_{2}^{2}
\end{array}\right)\otimes\left(\begin{array}{c}
d_{1}^{2}\\
d_{2}^{2}
\end{array}\right).
\end{eqnarray*}

To obtain the rank-2 canonical approximation, we minimize the functional
$\mathcal{F}$ 
\begin{eqnarray*}
\mathcal{F}=\frac{1}{2}\left\Vert \chi-\left(\begin{array}{c}
a_{1}^{1}\\
a_{2}^{1}
\end{array}\right)\otimes\left(\begin{array}{c}
b_{1}^{1}\\
b_{2}^{1}
\end{array}\right)\otimes\left(\begin{array}{c}
c_{1}^{1}\\
c_{2}^{1}
\end{array}\right)\otimes\left(\begin{array}{c}
d_{1}^{1}\\
d_{2}^{1}
\end{array}\right)+\left(\begin{array}{c}
a_{1}^{2}\\
a_{2}^{2}
\end{array}\right)\otimes\left(\begin{array}{c}
b_{1}^{2}\\
b_{2}^{2}
\end{array}\right)\otimes\left(\begin{array}{c}
c_{1}^{2}\\
c_{2}^{2}
\end{array}\right)\otimes\left(\begin{array}{c}
d_{1}^{2}\\
d_{2}^{2}
\end{array}\right)\right\Vert _{F}^{2}\\
=\frac{1}{2}\,\left(\left(\tau_{1}-\!(a_{1}^{1}b_{1}^{1}c_{1}^{1}d_{1}^{1}\!+\!a_{1}^{2}b_{1}^{2}c_{1}^{2}d_{1}^{2})\right)^{2}
+\!\left(\tau_{2}-\!(a_{2}^{1}b_{1}^{1}c_{1}^{1}d_{1}^{1}\!+\!a_{2}^{2}b_{1}^{2}c_{1}^{2}d_{1}^{2}\right)^{2}+\ldots
+\!\left(\tau_{16}-\!(a_{2}^{1}b_{2}^{1}c_{2}^{1}d_{2}^{1}\!+\!a_{2}^{2}b_{2}^{2}c_{2}^{2}d_{2}^{2}\right)^{2}\right).
\end{eqnarray*}

Let us denote 
\begin{eqnarray*}
A=\left[\begin{array}{cc}
a_{1}^{1} & a_{1}^{2}\\
a_{2}^{1} & a_{2}^{2}
\end{array}\right],\,B=\left[\begin{array}{cc}
b_{1}^{1} & b_{1}^{2}\\
b_{2}^{1} & b_{2}^{2}
\end{array}\right],\,C=\left[\begin{array}{cc}
c_{1}^{1} & c_{1}^{2}\\
c_{2}^{1} & c_{2}^{2}
\end{array}\right]\,\,\textrm{and}\,\,D=\left[\begin{array}{cc}
d_{1}^{1} & d_{1}^{2}\\
d_{2}^{1} & d_{2}^{2}
\end{array}\right].
\end{eqnarray*}

By ALS, $\mathcal{F}$ is minimized in an alternating way.
ALS first fixes $B,\,C$ and $D$ to minimize for $A,$ then
fixes $A,\,C$ and $D$ to minimize for $B$, then fixes $A,\,B$ and
$D$ to minimize for $C$ and finally fixes $A,\,B$ and $C$ to minimize
for $D.$ Since we are fixing all but one direction in each step
of an iteration, the problem reduces to a linear least-squares problem.
All the steps of one iteration are described below.

{\bf Step 1:} Fix $B,C$ and $D$ and solve for $A.$ The minimization leads
to the equations 
\begin{eqnarray*}
\frac{\partial\mathcal{F}}{\partial a_{1}^{1}}=0,\frac{\partial\mathcal{F}}{\partial a_{1}^{2}}=0\quad\text{and}
\quad\frac{\partial\mathcal{F}}{\partial a_{2}^{1}}=0,\frac{\partial\mathcal{F}}{\partial a_{2}^{2}}=0,
\end{eqnarray*}

which give a decoupled diagonal system 
\begin{eqnarray}
\left[\begin{array}{cc}
\hat{A}^{T}\hat{A} & 0\\
0 & \hat{A}^{T}\hat{A}
\end{array}\right]\left[\begin{array}{c}
\left[\begin{array}{c}
a_{1}^{1}\\
a_{1}^{2}
\end{array}\right]\\
\left[\begin{array}{c}
a_{2}^{1}\\
a_{2}^{2}
\end{array}\right]
\end{array}\right]=\left[\begin{array}{c}
\hat{A}^{T}\chi(1,:,:,:)\\
\hat{A}^{T}\chi(2,:,:,:)
\end{array}\right], 
\end{eqnarray}

where (see A2)

\begin{eqnarray}
\hat{A}=\left[\begin{array}{cc}
d_{1}^{1}c_{1}^{1}b_{1}^{1} & \,\,\,d_{1}^{2}c_{1}^{2}b_{1}^{2}\\
d_{1}^{1}c_{1}^{1}b_{2}^{1} & \,\,\,d_{1}^{2}c_{1}^{2}b_{2}^{2}\\
d_{1}^{1}c_{2}^{1}b_{1}^{1} & \,\,\,d_{1}^{2}c_{2}^{2}b_{1}^{2}\\
d_{1}^{1}c_{2}^{1}b_{2}^{1} & \,\,\,d_{1}^{2}c_{2}^{2}b_{2}^{2}\\
d_{2}^{1}c_{1}^{1}b_{1}^{1} & \,\,\,d_{2}^{2}c_{1}^{2}b_{1}^{2}\\
d_{2}^{1}c_{1}^{1}b_{2}^{1} & \,\,\,d_{2}^{2}c_{1}^{2}b_{2}^{2}\\
d_{2}^{1}c_{2}^{1}b_{1}^{1} & \,\,\,d_{2}^{2}c_{2}^{2}b_{1}^{2}\\
d_{2}^{1}c_{2}^{1}b_{2}^{1} & \,\,\,d_{2}^{2}c_{2}^{2}b_{2}^{2}
\end{array}\right]=D\odot C\odot B,\quad\chi(1,:,:,:)=\left[\begin{array}{c}
\tau_{1}\\
\tau_{3}\\
\tau_{5}\\
\tau_{7}\\
\tau_{9}\\
\tau_{11}\\
\tau_{13}\\
\tau_{15}
\end{array}\right],\quad\chi(2,:,:,:)=\left[\begin{array}{c}
\tau_{2}\\
\tau_{4}\\
\tau_{6}\\
\tau_{8}\\
\tau_{10}\\
\tau_{12}\\
\tau_{14}\\
\tau_{16}
\end{array}\right]. 
\end{eqnarray}

{\bf Step 2:} Fix $A,C,D$ and solve for $B.$ Then the equations,

\begin{eqnarray*}
\frac{\partial\mathcal{F}}{\partial b_{1}^{1}}=0,\frac{\partial\mathcal{F}}{\partial b_{1}^{2}}=0\quad\textrm{and}
\quad\frac{\partial\mathcal{F}}{\partial b_{2}^{1}}=0,\frac{\partial\mathcal{F}}{\partial b_{2}^{2}}=0
\end{eqnarray*}
give the linear system

\begin{eqnarray*}
\left[\begin{array}{cc}
\hat{B}^{T}\hat{B} & 0\\
0 & \hat{B}^{T}\hat{B}
\end{array}\right]\left[\begin{array}{c}
\left[\begin{array}{c}
b_{1}^{1}\\
b_{1}^{2}
\end{array}\right]\\
\left[\begin{array}{c}
b_{2}^{1}\\
b_{2}^{2}
\end{array}\right]
\end{array}\right]=\left[\begin{array}{c}
\hat{B}^{T}\chi(:,1,:,:)\\
\hat{B}^{T}\chi(:,2,:,:)
\end{array}\right]
\end{eqnarray*}

with 
\begin{eqnarray}
\hat{B}=D\odot C\odot A\quad\text{and}\quad\chi(:,1,:,:)=\left[\begin{array}{c}
\tau_{1}\\
\tau_{2}\\
\tau_{5}\\
\tau_{6}\\
\tau_{9}\\
\tau_{10}\\
\tau_{13}\\
\tau_{14}
\end{array}\right],\quad\chi(:,2,:,:)=\left[\begin{array}{c}
\tau_{3}\\
\tau_{4}\\
\tau_{7}\\
\tau_{8}\\
\tau_{11}\\
\tau_{12}\\
\tau_{15}\\
\tau_{16}
\end{array}\right].
\end{eqnarray}

{\bf Step 3:} Fix $A,B,D$ and solve for $C.$ Then the equations

\begin{eqnarray*}
\frac{\partial\mathcal{F}}{\partial c_{1}^{1}}=0,\frac{\partial\mathcal{F}}{\partial c_{1}^{2}}=0\quad\textrm{and}
\quad\frac{\partial\mathcal{F}}{\partial c_{2}^{1}}=0,\frac{\partial\mathcal{F}}{\partial c_{2}^{2}}=0
\end{eqnarray*}
give the linear system

\begin{eqnarray*}
\left[\begin{array}{cc}
\hat{C}^{T}\hat{C} & 0\\
0 & \hat{C}^{T}\hat{C}
\end{array}\right]\left[\begin{array}{c}
\left[\begin{array}{c}
c_{1}^{1}\\
c_{1}^{2}
\end{array}\right]\\
\left[\begin{array}{c}
c_{2}^{1}\\
c_{2}^{2}
\end{array}\right]
\end{array}\right]=\left[\begin{array}{c}
\hat{C}^{T}\chi(:,:,1,:)\\
\hat{C}^{T}\chi(:,:,2,:)
\end{array}\right]
\end{eqnarray*}

with 
\begin{eqnarray}
\hat{C}=D\odot B\odot A\quad\textrm{and}\quad\chi(:,:,1,:)=\left[\begin{array}{c}
\tau_{1}\\
\tau_{2}\\
\tau_{3}\\
\tau_{4}\\
\tau_{9}\\
\tau_{10}\\
\tau_{11}\\
\tau_{12}
\end{array}\right],\quad\chi(:,:,2,:)=\left[\begin{array}{c}
\tau_{5}\\
\tau_{6}\\
\tau_{7}\\
\tau_{8}\\
\tau_{13}\\
\tau_{14}\\
\tau_{15}\\
\tau_{16}
\end{array}\right]. 
\end{eqnarray}

{\bf Step 4:} Fix $A,B,C$ and solve for $D.$ The equations

\begin{eqnarray*}
\frac{\partial\mathcal{F}}{\partial d_{1}^{1}}=0,\frac{\partial\mathcal{F}}{\partial d_{1}^{2}}=0\quad\textrm{and}
\quad\frac{\partial\mathcal{F}}{\partial d_{2}^{1}}=0,\frac{\partial\mathcal{F}}{\partial d_{2}^{2}}=0
\end{eqnarray*}
give the linear system

\begin{eqnarray*}
\left[\begin{array}{cc}
\hat{D}^{T}\hat{D} & 0\\
0 & \hat{D}^{T}\hat{D}
\end{array}\right]\left[\begin{array}{c}
\left[\begin{array}{c}
d_{1}^{1}\\
d_{1}^{2}
\end{array}\right]\\
\left[\begin{array}{c}
d_{2}^{1}\\
d_{2}^{2}
\end{array}\right]
\end{array}\right]=\left[\begin{array}{c}
\hat{D}^{T}\chi(:,:,:,1)\\
\hat{D}^{T}\chi(:,:,:,2)
\end{array}\right]
\end{eqnarray*}

with 
\begin{eqnarray}
\hat{D}=C\odot B\odot A\quad\textrm{and}\quad\chi(:,:,:,1)=\left[\begin{array}{c}
\tau_{1}\\
\tau_{2}\\
\tau_{3}\\
\tau_{4}\\
\tau_{5}\\
\tau_{6}\\
\tau_{7}\\
\tau_{8}
\end{array}\right],\quad\chi(:,:,:,2)=\left[\begin{array}{c}
\tau_{9}\\
\tau_{10}\\
\tau_{11}\\
\tau_{12}\\
\tau_{13}\\
\tau_{14}\\
\tau_{15}\\
\tau_{16}
\end{array}\right]. 
\end{eqnarray}

Here the matrices $\hat{A,}\hat{B,}\hat{C}$ and $\hat{D}$ are of
size $8\times2.$ This is $2^{4-1}\times r$ where $r=2.$ But the matrices like $\hat{A}^{T}A$ appearing in the decoupled linear
systems are of very small size $2\times2$ for $r=2.$ Also one can
see that the matrices $\hat{A}^{T}\hat{A},\,\hat{B}^{T}\hat{B},\,\hat{C}^{T}\hat{C},\,\hat{D}^{T}\hat{D}$
are symmetric and positive definite.
\end{document}